\documentclass[12pt]{article}
\usepackage{amssymb}
\usepackage{amsmath}
\usepackage{latexsym}
\usepackage{mathrsfs}

\numberwithin{equation}{section}
\newtheorem{theorem}[equation]{Theorem}

\newtheorem{lemma}[equation]{Lemma}
\newtheorem{corollary}[equation]{Corollary}

\title{On the eigenvalues of a biharmonic Steklov problem}

\author{Davide Buoso\,\,   and Luigi Provenzano\footnote{Corresponding author}}

\date{\ }

\begin{document}

\newcommand{\rea}{\mathbb{R}}

\maketitle

%
%
%

\noindent
{\bf Abstract:} We consider an eigenvalue problem for the biharmonic operator with Steklov-type boundary conditions. We obtain it as a limiting Neumann problem for the biharmonic operator in a process of mass concentration at the boundary. We study the dependence of the spectrum upon the domain. We show analyticity of the symmetric functions of the eigenvalues under isovolumetric perturbations and prove that balls are critical points for such functions under measure constraint. Moreover, we show that the ball is a maximizer for the first positive eigenvalue among those domains with a prescribed fixed measure.
\vspace{11pt}

\vspace{6pt}


\newcommand{\oH}{{\mathaccent'27 H}}

\section{Introduction}
Let $\Omega$ be a bounded domain (i.e., a bounded connected open set) of class $C^2$ in $\mathbb R^N$, $N\geq 2$ and $\tau>0$. We consider the following Steklov eigenvalue problem for the biharmonic operator
\begin{eqnarray}\label{steklov}\left\{
\begin{array}{ll}
\Delta^2u-\tau\Delta u=0,& \ \ \ {\rm in}\ \Omega\,,\\
\frac{\partial^2 u}{\partial\nu^2}=0,& \ \ \ {\rm on}\ \partial\Omega\,,\\
\tau\frac{\partial u}{\partial \nu}-{\rm div}_{\partial\Omega}\left(D^2u.\nu\right)-\frac{\partial\Delta u}{\partial\nu}=\lambda u,&\ \ \ {\rm on}\ \partial\Omega\,,
\end{array}
\right.
\end{eqnarray}
in the unknowns $\lambda$ (the eigenvalue) and $u$ (the eigenfunction). Here $\nu$ denotes the unit outer normal to $\partial\Omega$, ${\rm div}_{\partial\Omega}$ the tangential divergence operator and $D^2u$ the Hessian matrix of $u$. The spectrum consists of a diverging sequence of eigenvalues of finite multiplicity
$$
0=\lambda_1<\lambda_2\leq\cdots\leq\lambda_j\leq\cdots,
$$
where we agree to repeat the eigenvalues according to their multiplicity.

When $N=2$, problem (\ref{steklov}) arises in the study of the vibration modes of a free elastic plate subjet to lateral tension (represented by the parameter $\tau$) whose total mass is concentrated at the boundary. We can describe this concentration phenomenon as follows.

For any $\varepsilon$ sufficiently small we consider the $\varepsilon$-neighborhood of $\partial\Omega$, namely $\omega_{\varepsilon}=\left\{x\in\Omega\,:\,0<d(x,\partial\Omega)<\varepsilon\right\}$. We fix $M>0$ and define the function $\rho_{\varepsilon}$ on $\Omega$ as follows

\begin{eqnarray*}\rho_{\varepsilon=}\left\{
\begin{array}{ll}
\varepsilon ,& \ \ \ {\rm in}\ \Omega\setminus\overline{\omega}_{\varepsilon}\,,\\
\frac{M-\varepsilon|\Omega\setminus\overline{\omega}_{\varepsilon}|}{|\omega_{\varepsilon}|},&\ \ \ {\rm in}\ \omega_{\varepsilon}\,.
\end{array}
\right.
\end{eqnarray*}
For any $x\in\Omega$ we have $\rho_{\varepsilon}(x)\rightarrow 0$ as $\varepsilon\rightarrow 0$. Moreover, $\int_{\Omega}\rho_{\varepsilon}=M$ for all $\varepsilon>0$. Then we consider the following eigenvalue problem for the biharmonic operator subject to Neumann boundary conditions
\begin{eqnarray}\label{neumann}\left\{
\begin{array}{ll}
\Delta^2u-\tau\Delta u={\lambda(\varepsilon)}\rho_{\varepsilon} u ,& \ \ \ {\rm in}\ \Omega\,,\\
\frac{\partial^2 u}{\partial\nu^2}=0,& \ \ \ {\rm on}\ \partial\Omega\,,\\
\tau\frac{\partial u}{\partial\nu}-{\rm div}_{\partial\Omega}\left(D^2u.\nu\right)-\frac{\partial\Delta u}{\partial\nu}=0,&\ \ \ {\rm on}\ \partial\Omega\,.
\end{array}
\right.
\end{eqnarray}
The spectrum consists of a diverging sequence of eigenvalues of finite multiplicity
$$
0=\lambda_1(\varepsilon)<\lambda_2(\varepsilon)\leq\cdots\leq\lambda_j(\varepsilon)\leq\cdots,
$$
where we agree to repeat the eigenvalues according to their multiplicity. Here we emphasize the dependence of the eigenvalues on the parameter $\varepsilon$.

We remark that for $N=2$ problem (\ref{neumann}) provides the fundamental modes of vibration of a free elastic plate with mass density $\rho_{\varepsilon}$ and total mass $M$, as discussed in \cite[Chasman]{chas1}. We refer to \cite{chas1} for the derivation and the physical interpretation of problem (\ref{neumann}). 

It is possible to prove that the eigenvalues and the eigenfunctions of (\ref{neumann}) converge to the eigenvalues and eigenfunctions of (\ref{steklov}) as $\varepsilon$ goes to zero (see, e.g., \cite{arr,buproz,proz}).

The aim of this paper is to study a few properties concerning the dependence of the eigenvalues of (\ref{steklov}) upon perturbations of the domain $\Omega$ which preserve the measure.

First, we study the asymptotic behavior of the eigenvalues of (\ref{neumann}) as $\varepsilon\rightarrow 0$ providing an interpretation of (\ref{steklov}) as the model of a free vibrating plate with all the mass concentrated at the boundary (see Theorem \ref{eigenconvergence}). This fact suggests that (\ref{steklov}) is the natural fourth order generalization of the classical Steklov eigenvalue problem for the Laplace operator, see \cite{stek} and the recent \cite{la14} for related problems.

Second, we consider the problem of the optimal shape of $\Omega$ for the eigenvalues of (\ref{steklov}) under the constraint that the measure of $\Omega$ is fixed. This problem has been largely investigated for the Laplace operator subject to different homogeneous boundary conditions. We refer to \cite[Henrot]{he} for a collection of results on the subject. See also \cite[Bandle]{bandle}. As far as the biharmonic operator is concerned, only a few results exist in literature. It has been proved in  \cite[Nadirashvili]{nadir} for $N=2$ and soon generalized in \cite[Ashbaugh, Benguria]{ash} for $N=3$ that the ball is a minimizer for the first eigenvalue of the biharmonic operator subject to Dirichlet boundary conditions. In the recent paper \cite{chas1}, it has been proved that the first positive eigenvalue of problem (\ref{neumann}) with constant mass density $\rho\equiv 1$ is maximized by the ball among those sets with a fixed measure.\\
As for Steklov boundary conditions, we refer to \cite[Bucur, Ferrero, Gazzola]{bucurgazzola} and the references therein. The authors consider the following eigenvalue problem
\begin{eqnarray}\label{steklovgaz}\left\{
\begin{array}{ll}
\Delta^2u=0,& \ \ \ {\rm in}\ \Omega\,,\\
u=0,& \ \ \ {\rm on}\ \partial\Omega\,,\\
\Delta u=\lambda \frac{\partial u}{\partial\nu},&\ \ \ {\rm on}\ \partial\Omega\,.
\end{array}
\right.
\end{eqnarray}
Problem (\ref{steklovgaz}) should not be confused with problem (\ref{steklov}) and reveals a rather different nature. (We note that one may refer to Steklov-type boundary conditions for those problems where the spectral parameter enters the boundary conditions.) 

By following the approach developed in \cite{bula2013,buosohinged} we prove that simple eigenvalues and the symmetric functions of multiple eigenvalues of (\ref{steklov}) depend real analytically upon transformations of the domain $\Omega$ (see Theorem \ref{symmetric}) and we characterize those critical trasformations which preserve the measure (see Corollary \ref{crit}). See also \cite{lala2004, lalacri, lala2007}. Then we show that the ball is a critical point for all simple eigenvalues and all symmetric functions of the eigenvalues under measure constraint in the sense of Theorem \ref{ballcrit}.

Finally, we prove the following isoperimetric inequality: {\it ``The ball is a maximizer for the first positive eigenvalue of problem (\ref{steklov}) among those bounded domains with a fixed measure''} (see Theorem \ref{isotheorem}). To do so, we follow the approach of \cite{chas1} and in particular we study problem (\ref{steklov}) when $\Omega$ is the unit ball in $\mathbb R^N$, identifying the first positive eigenvalue and the corresponding eigenfunctions.

Detailed proofs of the results announced in this paper can be found in \cite{buproz}.

\section{Asymptotic behavior of Neumann eigenvalues}

Let $\Omega$ be a bounded domain in $\mathbb R^N$ of class $C^2$. Let $\lambda_j$ and $\lambda_j(\varepsilon)$, $j\in\mathbb N\setminus\left\{0\right\}$, be the eigenvalues of (\ref{steklov}) and (\ref{neumann}) respectively. For the sake of simplicity and without any loss of generality we assume that $M=|\partial\Omega|$. We recall that $\lambda_1=\lambda_1(\varepsilon)=0$, while $\lambda_2,\lambda_2(\varepsilon)>0$ for all $\varepsilon>0$. 

We have the following result concerning the spectral convergence of problem (\ref{neumann}) to problem (\ref{steklov}).
\begin{theorem}\label{eigenconvergence}
Let $\Omega$ be a bounded domain in $\mathbb R^N$ of class $C^2$. Then $\lambda_j(\varepsilon)\rightarrow\lambda_j$ for all $j\in\mathbb N\setminus\left\{0\right\}$. Moreover the projections on the eigenspaces associated with the eigenvalues converge in norm.
\end{theorem}

This theorem can be proved by using the notion of compact convergence for the resolvent operators which implies, in the case of selfadjoint operators, convergence in norm. It is well known that if a family of selfadjoint operators $A_{\varepsilon}$ converges in norm to a selfadjoint operator $A$, then isolated eigenvalues of $A$ are exactly the limits of eigenvalues of $A_{\varepsilon}$ counting multiplicity. Moreover, eigenprojections converge in norm. We refer to \cite{buproz,proz} for more details. We also refer to the recent paper \cite{ArLa} for a general approach to the shape sensitivity analysis of higher order operators.

Theorem \ref{eigenconvergence} justifies our interpretation of problem (\ref{steklov}) as the equations of a free vibrating plate whose mass is concentrated at the boundary. However, we can also directly obtain (\ref{steklov}) by deriving the equations of motions of a free vibrating plate with constant surface density. To do so, we follow the approach of \cite[ch.10-8]{wein} in the case $N=2$. We represent the displacement at rest of the plate by means of a domain $\Omega\subset\mathbb R^2$ and we describe the vertical deviation from the equilibrium during the vibration of each point $(x,y)\in\Omega$ at time $t$ by means of a function $v(x,y,t)\in C^2(\Omega\times[t_1,t_2])$. Then we write the Hamiltonian $\mathcal H$ of the system
\begin{equation}\label{hamilton}
\mathcal H=\frac{1}{2}\int_{t_1}^{t_2}\int_{\partial\Omega}\dot v^2d\sigma dt-\frac{1}{2}\int_{t_1}^{t_2}\int_{\Omega}\left(v^2_{xx}+v^2_{yy}+2v^2_{xy}\right)+\tau\left(v_x^2+v_y^2\right)dxdydt.
\end{equation}
According to Hamilton's Variational Principle, we have to minimize the Hamiltonian $\mathcal H$. Let $v\in C^2(\Omega\times[t_1,t_2])$ be a minimizer for $\mathcal H$. Then by differentiating (\ref{hamilton}) it follows that $v$ satisfies
\begin{eqnarray*}
&&-\int_{t_1}^{t_2}\int_{\partial\Omega}\eta\ddot v d\sigma dt-\int_{t_1}^{t_2}\int_{\Omega}\eta\left(\Delta^2v-\tau\Delta v\right)dxdydt\\
&-&\int_{t_1}^{t_2}\int_{\partial\Omega}\frac{\partial\eta}{\partial\nu}\frac{\partial^2v}{\partial\nu^2}-\eta\left(\tau\frac{\partial v}{\partial\nu}-{\rm div}_{\partial\Omega}\left(D^2v.\nu\right)_{\partial\Omega}-\frac{\partial\Delta v}{\partial\nu}\right)d\sigma dt=0\nonumber,
\end{eqnarray*}
for all $\eta\in C^2(\Omega\times[t_1,t_2])$. We refer \cite{chas1} for the details. By the arbitrary choice of $\eta$ we obtain
\begin{equation*}
   \begin{cases}
	\Delta^2v-\tau\Delta v=0, & {\rm in}\ \Omega,\\
   \frac{\partial^2v}{\partial\nu^2}=0, & {\rm on}\ \partial\Omega,\\
	\ddot v+\tau\frac{\partial v}{\partial\nu}-{\rm div}_{\partial\Omega}\left(D^2v.\nu\right)-\frac{\partial\Delta v}{\partial\nu}=0, & {\rm on}\ \partial\Omega,
   \end{cases}
\end{equation*}
for all $t\in\mathbb R$. As is customary, by looking for solution of the form $v(x,y,t)=u(x,y)\psi(t)$. We find that the temporal component $\psi(t)$ solves the ordinary differential equation $-\ddot\psi=\lambda\psi$ for all $t\in[t_1,t_2]$, while the spatial component $u$ solves problem (\ref{steklov}).

\section{Isovolumetric perturbations}
Given a bounded domain in $\mathbb R^N$ of class $C^2$, we set

	$$\Phi(\Omega)=\left\{\phi\in\left(C^2\left(\overline{\Omega}\right)\right)^N:\phi\ {\rm injective}\ {\rm and}\ \inf_{\Omega}|\det D\phi|>0\right\}.$$
We observe that if $\Omega$ is of class $C^2$ and $\phi\in\Phi(\Omega)$, it makes sense to study problem (\ref{steklov}) on $\phi(\Omega)$. For any $\phi\in\Phi(\Omega)$ we denote by $\lambda_j(\phi)$, $j\in\mathbb N\setminus\left\{0\right\}$, the eigenvalues of (\ref{steklov}) on $\phi(\Omega)$.

We plan to study the dependence of the eigenvalues upon the function $\phi$. In general, one cannot expect differentiability of the eigenvalues with respect to $\phi$. This is due, for example, to well known bifurcation phenomena that occur when multiple eigenvalues split from a simple eigenvalue. However, as is pointed out in \cite{bula2013,buosohinged}, in the case of multiple eigenvalues it is possibile to prove analyticity for the symmetric functions of the eigenvalues. Namely, given a finite set of indexes $F\subset\mathbb N\setminus\left\{0\right\}$, one can consider the symmetric functions of the eigenvalues with indexes in $F$
\begin{equation*}
\Lambda_{F,s}(\phi)=\sum_{j_1<\cdots<j_s\in F}\lambda_{j_1}(\phi)\cdots\lambda_{j_s}(\phi),
\end{equation*}
and prove that such functions are real analytic on the set
\begin{equation*}
\mathcal A_{\Omega}[F]=\left\{\phi\in\Phi(\Omega)\,:\,\lambda_l(\phi)\not\in\left\{\lambda_j(\phi):j\in F\right\}\,\forall l\in\mathbb N\setminus\left(F\cup\left\{0\right\}\right)\right\}.
\end{equation*}
Then it is possible to find formulas for the Fr\'echet derivatives of the symmetric functions of the eigenvalues. It is convenient to set
\begin{equation*}
\Theta_{\Omega}[F]=\left\{\phi\in\mathcal A_{\Omega}[F]\,:\,\lambda_{j_1}(\phi)=\lambda_{j_2}(\phi)\,,\forall j_1,j_2\in F\right\}.
\end{equation*}
\begin{theorem}\label{symmetric}
Let $\Omega$ be a bounded domain of $\mathbb R^N$ of class $C^2$. Let $F$ be a finite non-empty subset of $\mathbb N\setminus\left\{0\right\}$. Then $\mathcal A_{\Omega}$ is open in $\Phi(\Omega)$ and $\Lambda_{F,s}$ are real analytic in $\mathcal A_{\Omega}$. Moreover, let $\tilde\phi\in\Theta_{\Omega}[F]$ be such that $\partial\tilde\phi(\Omega)\in C^4$. Let $v_1,...,v_{|F|}$ be a orthonormal basis of the eigenspace associated with the eigenvalue $\lambda_F(\tilde\phi)$. Then
\begin{multline}
		\label{deriv}
			d|_{\phi=\tilde{\phi}}(\Lambda_{F,s})(\psi)=-\lambda_F^s(\tilde{\phi})\binom{|F|-1}{s-1}
\sum_{l=1}^{|F|}\int_{\partial\tilde{\phi}(\Omega)}\Big(\lambda_FKv_l^2\\
+\lambda_F\frac{\partial(v_l^2)}{\partial\nu}-\tau|\nabla v_l|^2-|D^2v_l|^2\Big)\mu\cdot\nu d\sigma,
		\end{multline}
	for all $\psi\in(C^2(\Omega))^N$, where $\mu=\psi\circ\phi^{(-1)}$, and $K$ denotes the mean curvature on $\partial\tilde{\phi}(\Omega)$.
\end{theorem}
The proof follows the lines of the corresponding results provided in \cite{bula2013} and \cite{buosohinged} for general poly-harmonic operators subject to Dirichlet boundary conditions and for the biharmonic operator subject to hinged boundary conditions.

We consider now the problem of finding critical points $\phi\in\Phi(\Omega)$ for the symmetric functions of the eigenvalues under the condition that $\phi$ preserves the measure. We set $\mathcal V(\phi)=\int_{\phi(\Omega)}dy=\int_{\Omega}|{\rm det D\phi}|dx$. We fix $\mathcal V_0\in]0,+\infty[$ and consider the set $V(\mathcal V_0)=\left\{\phi\in\Phi(\Omega)\,:\,\mathcal V(\phi)=\mathcal V_0\right\}$. Given $\Omega$ such that $|\Omega|=\mathcal V_0$, $V(\mathcal V_0)$ is the subset of $\Phi(\Omega)$ of those functions $\phi$ preserving the measure. By formula (\ref{deriv}) and by the Lagrange Multipliers Theorem we can characterize the critical points.
\begin{corollary}\label{crit}
Let all the assumptions of Theorem \ref{symmetric} hold. Then $\tilde\phi\in\Phi(\Omega)$ is a critical point for $\Lambda_{F,s}$ on $V(\mathcal V_0)$ if and only if there exists a constant $c\in\mathbb R$ such that
\begin{equation*}
	\sum_{l=1}^{|F|}\left(\lambda_F(\tilde{\phi})\left(Kv_l^2+\frac{\partial v_l^2}{\partial\nu}\right)-\tau|\nabla v_l|^2-|D^2v_l|^2\right)=c,\,\,{\rm a.e.}\,{\rm on}\ \partial\tilde\phi(\Omega),
\end{equation*}
where $K$ denotes the mean curvature on $\partial\tilde\phi (\Omega)$.
\end{corollary}
Thanks to Corollary \ref{crit} we can prove that balls are critical points for the symmetric functions of the eigenvalues under measure constraint, in the sense of the following
\begin{theorem}\label{ballcrit}
 Let $\tilde{\phi}\in\Phi(\Omega)$ be such that
	$\tilde{\phi}(\Omega)$ is a ball. Let $\tilde{\lambda}$ be an eigenvalue of the problem in $\tilde{\phi}(\Omega)$,
	and let $F$ be the set of all $j\in\mathbb{N}\setminus\{0\}$ such that $\lambda_j(\tilde{\phi})=\tilde{\lambda}$.
	Then $\Lambda_{F,s}$ has a critical point at $\tilde{\phi}$ on $V(\mathcal{V}_0)$,
	for all $s=1,\dots,|F|$.
\end{theorem}
The proof can be carried out as in \cite{bula2013,buosohinged}. Namely, given $\lambda$ an eigenvalue of problem (\ref{steklov}) on the unit ball $B$ in $\mathbb R^N$, consider the subset $F$ of $\mathbb N\setminus\left\{0\right\}$ of those indexes $j$ such that the $j$-th eigenvalue of problem (\ref{steklov}) in $B$ coincides with $\lambda$. Consider then $v_1,...,v_{|F|}$ an orthonormal basis of the eigenspace associated with the eigenvalue $\lambda$, where the orthonormality is taken with respect to the scalar product in $L^2(\partial B)$. Then it is possible to show that the quantities $\sum_{j=1}^{|F|}v_j^2$, $\sum_{j=1}^{|F|}|\nabla v_j|^2$ and $\sum_{j=1}^{|F|}|D^2v_j|^2$ are radial functions. This fact and the fact that the mean curvature $K$ is constant on the ball allow to conclude.

\section{The isoperimetric inequality}

Let us consider problem (\ref{steklov}) when $\Omega=B$ is the unit ball in $\mathbb R^N$. It is convenient to use spherical coordinates $\left(r,\theta\right)$ in $\mathbb R^N$, where $\theta=\left(\theta_1,...,\theta_{N-1}\right)$, with $r\in \left[0,1\right[$\,, ${\theta_1},...,{\theta_{N-2}}\in \left[0,\pi\right]$\,,  $\theta_{N-1}\in\left[0,2\pi\right]$. In this case the boundary conditions can we written in the following form

\begin{equation}\label{Steklov-Bi-Ball}
\left\{\begin{array}{ll}
\frac{\partial^2 u}{\partial r^2 }_{|_{r=1}}=0 ,\\
\tau\frac{\partial u}{\partial r }-\frac{1}{r^2}{\Delta_S}\Big(\frac{\partial u}{\partial r}-\frac{u}{r}\Big)-\frac{\partial\Delta u}{\partial r}_{|_{r=1}}=\lambda u_{|_{r=1}},
\end{array}\right.
\end{equation}
where $\Delta_S$ is the angular part of the Laplacian (see \cite{chas1} for details). Then, the eigenfunctions of problem (\ref{steklov}) on the ball can be described explicitly as in the following lemma. We refer to \cite[ch.9]{abram} for well-known definitions and properties of Bessel functions. 

\begin{lemma}\label{eigenfunctions}
Let $B$ be the unit ball in $\mathbb R^N$. An eigenfunction $u$ of (\ref{steklov}) is of the form $u(r,\theta)=R_l(r)Y_l(\theta)$, where $Y_l(\theta)$ is a spherical harmonic of some order $l\in\mathbb N$,
\begin{equation}\label{index}
R_l(r)=A_lr^l+B_li_l(\sqrt{\tau}r)
\end{equation}
and $A_l$ and $B_l$ are suitable constants such that
\begin{equation}\label{combination}
B_l=\frac{l(1-l)}{\tau i_l''(\sqrt{\tau})}A_l.
\end{equation}
Here $i_l$ denotes the ultraspherical modified Bessel function of the first kind, which is defined by
$$
i_l(z)=z^{1-\frac{N}{2}}I_{\frac{N}{2}-1+l}(z),
$$
where $I_l(z)$ denotes the modified Bessel function of the first kind.
\end{lemma}

We note that equality (\ref{combination}) is obtained by imposing the boundary conditions (\ref{Steklov-Bi-Ball}) to the function (\ref{index}). 

Lemma \ref{eigenfunctions} allows to find explicit formulas for the eigenvalues. In the sequel we will denote by $\lambda_{(l)}$ the eigenvalue corresponding to the eigenfunction $u_l$ defined in Lemma \ref{eigenfunctions}.
 
\begin{lemma}\label{eigenlem}
The eigenvalues $\lambda_{(l)}$ of problem (\ref{steklov}) on $B$ are delivered by the formula
\small
\begin{eqnarray*}
\lambda_{(l)}&=&l\Big((1 - l) l i_l(\sqrt{\tau}) + 
 \tau i_l''(\sqrt{\tau})\Big)^{-1}\Big[ 3 (l-1) l (l+N-2) i_l(\sqrt{\tau})\nonumber\\
 &&- (l-1) \sqrt{
    \tau} \big(N-1 + 2 N l + 2l(l-2) l + \tau\big) i_l'(\sqrt{
     \tau})\nonumber\\
		&&+ \tau \big((l-1) (l+2N-3)+ \tau\big) i_l''(
     \sqrt{\tau})\nonumber\\
		&&+ (l-1) \tau\sqrt{\tau} 
i_l'''(\sqrt{\tau})\Big],
\end{eqnarray*}
\normalsize
with $l\in\mathbb N$.
\end{lemma}

Now we need to identify the index $l$ satisfying $\lambda_{(l)}=\lambda_2$, that is the first positive eigenvalue of (\ref{steklov}). This is done by means of the following
\begin{lemma}\label{fund}
The first positive eigenvalue of problem (\ref{steklov}) on $B$ is $\lambda_2=\lambda_{(1)}=\tau$. The corresponding eigenspace is generated by the coordinate functions $\left\{x_1,...,x_N\right\}$.
\end{lemma}
The proof of Lemma \ref{fund} consists in two steps. In the first step we observe that $0=\lambda_{(0)}<\lambda_{(1)}=\tau$. Moreover, by using well known recurrence relations for modified ultraspherical Bessel functions of the first kind and their derivatives we are able to prove that $\lambda_{(1)}<\lambda_{(2)}$. In the second step we show that for any smooth radial function $R(r)$, the Rayleigh quotient
$$
\mathcal Q(R(r)Y_l(\theta))=\frac{\int_{\Omega}|D^2(R(r)Y_l(\theta))|^2+\tau|\nabla(R(r)Y_l(\theta))|^2dx}{\int_{\partial \Omega}R(r)^2Y_l(\theta)^2d\sigma}
$$
is an increasing function of $l$ for $l\geq 2$. This, combined with the variational characterization of the eigenvalues, allows us to conclude that $\lambda_{(l)}$ is an increasing function of $l$ for $l\geq 2$.

We are ready to state the isoperimetric inequality.

\begin{theorem}\label{isotheorem}
Among all bounded domains of class $C^2$ with fixed measure, the ball maximizes the first non-negative eigenvalue, that is $\lambda_2(\Omega)\leq\lambda_2(\Omega^*)$, where $\Omega^*$ is a ball with the same measure as $\Omega$.
\end{theorem}

The proof can be carried out as in \cite[par.7.3]{he}. Namely, we use the following variational characterization of the sum of inverse of eigenvalues

\begin{equation}\label{varsum}
\sum_{l=2}^{N+1}\frac{1}{\lambda_l(\Omega)}=\max\Bigg\{\sum_{l=2}^{N+1}\int_{\partial\Omega}v_l^2d\sigma\Bigg\},
\end{equation}
where $\{v_l\}_{l=2}^{N+1}$ is a family in $H^2(\Omega)$ satisfying $\int_{\Omega}D^2v_i:D^2v_j+\tau\nabla v_i\cdot\nabla v_j dx=\delta_{ij}$ and $\int_{\partial\Omega}v_l d\sigma=0$ for all $l=2,...,N+1$. We plug the functions ${v_l}={(\tau|\Omega|)^{-\frac{1}{2}}x_l}$, with $l=1,...,N$, into (\ref{varsum}) and we use the following inequality 
\begin{equation}\label{isop}
\int_{\partial\Omega}f(|x|)d\sigma\geq\int_{\partial\Omega^*}f(|x|)d\sigma,
\end{equation}
where $\Omega^*$ is the ball with the same measure of $\Omega$ and $f$ is a continuous, non-negative, non-decreasing function defined on $[0,+\infty)$ and moreover is such that the map $t\mapsto\big(f(t^{1/N})-f(0)\big)t^{1-(1/N)}$ is convex. Then the  isoperimetric inequality easily follows. We refer to \cite{hile} for the proof of (\ref{varsum}) and to \cite{betta} for the proof of (\ref{isop}).

\section*{Acknowledgements}	
The authors are deeply thankful to Prof.\ Pier Domenico Lamberti who suggested the problem, and also for many useful discussions.
The authors acknowledge financial support from the research project
`Singular perturbation problems for differential operators' Progetto di Ateneo
of the University of Padova. The authors are members of the Gruppo Nazionale
per l'Analisi Matematica, la Probabilit\`a e le loro Applicazioni (GNAMPA) of
the Istituto Nazionale di Alta Matematica (INdAM).

\noindent {\small
Davide Buoso and Luigi Provenzano\\
Dipartimento di Matematica\\
Universit\`{a} degli Studi di Padova\\
Via Trieste,  63\\
35126 Padova\\
Italy\\
e-mail:	dbuoso@math.unipd.it \\
e-mail: proz@math.unipd.it

}
\end{document}